\documentclass{article}[12pt]
\usepackage{latexsym, graphics, amsfonts}
\begin{document}
\newtheorem{theorem}{Theorem}[section]
\newtheorem{corollary}[theorem]{Corollary}
\newtheorem{lemma}[theorem]{Lemma}
\newtheorem{proposition}[theorem]{Proposition}
\newcommand{\Match}[1]{\mathcal{M}(#1)}
\newcommand{\bx}[1]{\mathcal{B}(#1)}
\newcommand{\insertfig}[1]{\begin{center}\includegraphics{#1}\end{center}}
\newenvironment{proof}{{\bf Proof:}}{\hfill\rule{2mm}{2mm}}

\title{Applications of Graphical Condensation for
Enumerating Matchings and Tilings}

\author{Eric H. Kuo\footnote{Computer Science Division, University of
California, Berkeley, CA, 94720, USA.  E-mail:{\tt ekuo@cs.berkeley.edu.}
This material is based on work performed under the auspices of the
Undergraduate Research Opportunities Program at MIT.
The work was also supported by an NSF Graduate Fellowship.}}

\date{August 22, 2003}
\maketitle

\begin{abstract}
A technique called graphical condensation is used to prove various
combinatorial identities among numbers of (perfect) matchings of planar bipartite graphs and
tilings of regions.  Graphical condensation involves superimposing matchings
of a graph onto matchings of a smaller subgraph, and then re-partitioning the united
matching (actually a multigraph) into matchings of two other subgraphs,
in one of two possible ways.
This technique can be used to enumerate perfect matchings
of a wide variety of planar bipartite graphs.
Applications include domino tilings of Aztec diamonds and rectangles,
diabolo tilings of fortresses,
plane partitions, and transpose complement plane partitions.
\end{abstract}

\section{Introduction}\label{sec:intro}

\begin{figure}
\insertfig{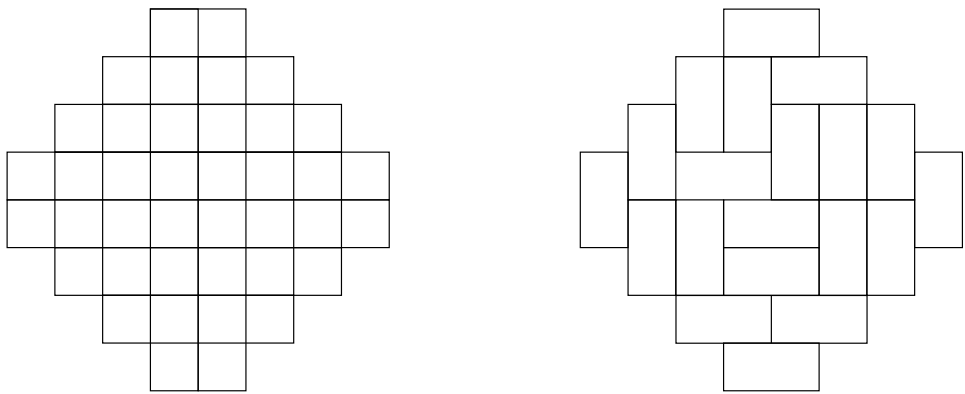}
\caption{Order 4 Aztec diamond and a sample domino tiling.}
\label{fig:4-ad}
\end{figure}

The {\it Aztec diamond\/} of order $n$ is defined as the union of all
unit squares whose corners are lattice points which lie within the
region $\{(x,y): |x| + |y| \leq n+1\}$.   A {\it domino \/} is simply a 1-by-2
or 2-by-1 rectangle whose corners are lattice points.  A {\it domino
tiling \/} of a region $R$ is a set of non-overlapping dominoes whose union
is $R$.  Figure~\ref{fig:4-ad} shows an Aztec diamond of order 4 and a sample
domino tiling.

In~\cite{GCZ}, it
was conjectured that the number of tilings for the order-$n$ Aztec diamond is
$2^{n(n+1)/2}$.  The conjecture was proved in~\cite{eklp}.
As the author went about trying to enumerate domino tilings for a similar
region, he discovered a new technique called {\it graphical
condensation.\/}  This technique has some far-reaching applications for
proving various combinatorial identities.  These identities usually take
the form
\[ T(R_1)T(R_2) = T(R_3)T(R_4) + T(R_5)T(R_6), \]
where $T(R_i)$ stands for the number of tilings for a region $R_i$.
In our applications, the regions $R_i$ are complexes
built out of vertices, edges, and faces,
and the legal tiles correspond to pairs of faces that share an edge;
a collection of such tiles constitutes a tiling if each face of $R_i$
belongs to exactly one tile in the collection
(such tilings are sometimes called diform tilings).
Each region $R_i$ could be represented by its dual graph $G_i$.  The number of
tilings for $R_i$ would equal the number of perfect matchings of $G_i$.
Thus we could replace each term $T(R_i)$ in the identity with
$M(G_i)$, which stands for the number of perfect matchings of $G_i$.
(Hereafter, it will be understood that any use of the term ``matching''
refers to a perfect matching.)

Graphical condensation involves superimposing a matching of one graph onto
a matching of another, and then partitioning that union into matchings of
two other graphs.  The phrase {\it graphical condensation} comes from the
striking resemblance between Dodgson condensation of determinants and
graphical condensation of Aztec diamonds.  A proof of Dodgson condensation
which illustrates this striking resemblance can be found in~\cite{zttmw}.

This article describes how graphical condensation can be used to prove
bilinear relations among numbers of matchings of planar bipartite
graphs or diform tilings of regions.
Among the applications are domino tilings of Aztec diamonds
(as well as some variant regions with holes in them),
and rhombus (or lozenge) tilings of semiregular hexagons
(equivalent to plane partitions),
with or without the requirement of bilateral symmetry.
The main result extends
to weighted enumeration of matchings of edge-weighted graphs,
and this extension gives us a simple way
to apply the method to count domino tilings of rectangles
and diabolo tilings of fortresses.

\section{Enumerative Relations Among matchings of
Planar Bipartite Graphs}\label{sec:pbg}

Before we state our enumerative relations, let us introduce some notation.
We will be working with a bipartite graph $G=(V_1,V_2,E)$ in which $V_1$ and
$V_2$ are disjoint sets of vertices in $G$ and every edge in $E$ connects a
vertex in $V_1$ to a vertex in $V_2$.  If $U$ is a subset of vertices in $G$,
then $G-U$ is the subgraph of $G$ obtained by deleting the vertices in $U$
and all edges incident to those vertices.  If $a$ is a vertex in $G$, then
 $G-a = G-\{a\}$.  Finally, we will let $M(G)$ be the number of perfect
matchings of $G$, and $\Match{G}$ be the set of all perfect matchings of $G$.

In order to state the enumerative relations, we must first embed $G$ into the
plane $\mathbb{R}^2$.  The {\it plane} graph $G$ divides $\mathbb{R}^2$ into
faces, one of which is unbounded.

\begin{figure}
\insertfig{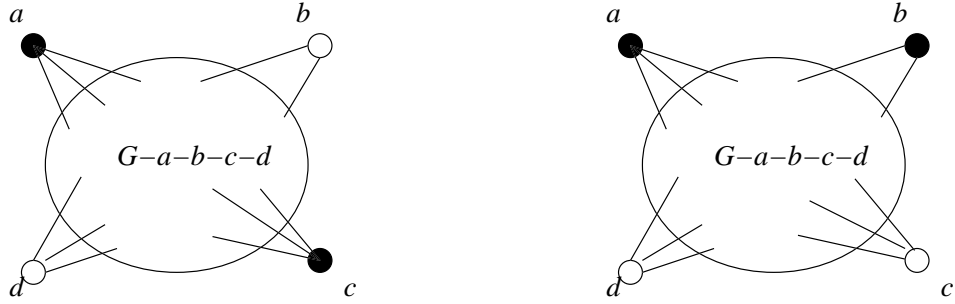}
\caption{Left: Graph for Theorem~\ref{thm:abcd}.
Right: Graph for Theorem ~\ref{thm:ad-bc}. }
\label{fig:ab-cd}
\end{figure}

\begin{theorem}\label{thm:abcd}
Let $G=(V_1,V_2,E)$ be a plane bipartite graph in which $|V_1|=|V_2|$.  Let
vertices $a$, $b$, $c$, and $d$ appear in a cyclic order on a face of $G$.  
(See Figure~\ref{fig:ab-cd}, left.  Note $a,b,c,d$ lie on the unbounded face.) 
If $a, c \in V_1$ and $b,d \in V_2$, then
\[ M(G) M(G-\{a,b,c,d\}) = M(G-\{a,b\})M(G-\{c,d\}) + M(G-\{a,d\})M(G-\{b,c\}). \]
\end{theorem}

\begin{proof}
To prove this relation, we would like to
establish that the two sets $\Match{G} \times \Match{G-\{a,b,c,d\}}$ and
$(\Match{G-\{a,b\}} \times \Match{G-\{c,d\}}) \cup (\Match{G-\{a,d\}} \times
\Match{G-\{b,c\}})$ have the same cardinality. Consider superimposing a
matching of $G-\{a,b,c,d\}$ onto a matching of $G$.  Whenever both matchings
share a common edge, we retain both edges and place a doubled edge in the
united matching.  Thus in the united matching (strictly speaking
a multigraph, since some edges may belong with multiplicity 2),
each vertex has degree 2 except for
$a, b, c,$ and $d$, which have degree 1.

Now consider superimposing a matching of $G-\{a,b\}$ onto a matching of
$G-\{c,d\}$.  Each vertex in the resulting graph has degree 2 except for $a,
b, c,$ and $d$, which have degree 1.  The same type of graph results from
superimposing a matching of $G-\{a,d\}$ onto a matching of $G-\{b,c\}$.

We define $\mathcal{H}$ to be the set of multigraphs on the vertices of $G$ in
which vertices $a,b,c,$ and $d$ have degree 1, and all remaining vertices 
have degree 2.  The edges of $G$ form cycles, doubled edges, and two paths 
whose endpoints are $a, b, c,$ and $d$.
Each pair of graphs in $\Match{G} \times \Match{G-\{a,b,c,d\}}$,
$\Match{G-\{a,b\}} \times \Match{G-\{c,d\}}$, and $\Match{G-\{a,d\}} \times
\Match{G-\{b,c\}}$ can be merged to form a multigraph in $\mathcal{H}$.
(Hereafter, we shall drop the prefix ``multi-'' and
refer to the elements of
$\mathcal{H}$ as simply graphs.)

Let $H$ be a graph in $\mathcal{H}$.  From $a$, we can trace a path through
$H$ until we hit another vertex of degree 1.  No vertex can be visited twice
by this path since each vertex has degree at most two.  Eventually we must
end at one of the other vertices of degree 1.  If one path connects $a$ to
$b$, then the path from $c$ must end at the remaining degree-1 vertex $d$.
Otherwise if $a$ connects to $d$, then $b$ must connect to $c$.  And since
$a, b, c,$ and $d$ occur in cyclic order around a face of $G$,
it is impossible for one path to connect $a$ to $c$ and the other path to
connect $b$ to $d$.  If such paths existed, then they would have to
intersect, forcing some other vertex to have a degree greater than 2.

We now show that $H$ can be partitioned into a matching $M_1$ of $G$
and a matching $M_2$ of $G-\{a,b,c,d\}$ in $2^k$ ways, where $k$ is the number
of cycles in $H$.  Since $H$ is bipartite, each cycle
has even length.  We partition each cycle in $H$ so that adjacent edges go
into different matchings; each vertex in a cycle is incident to one edge from
each matching.  Each doubled edge is split and shared between both matchings.
Since the paths connect $a$ to $b$ (or $d$) and $c$ to $d$ (or $b$),
one end of each path must belong to $V_1$ and the other end must be in
$V_2$.  Thus each path has odd length (as measured by the number of edges),
so we may assign the edges at the ends of each path to $M_1$.  The 
remaining edges in the paths are assigned to $M_1$ and $M_2$,
and thus it is always possible to partition $H$ into matchings $M_1$ and $M_2$.
Since there are two choices for distributing edges in each cycle of $H$ into
matchings $M_1$ and $M_2$, there are $2^k$ possible ways to partition $H$ into
matchings of $G$ and $G-\{a,b,c,d\}$.

Next, we show that $H$ can always be partitioned into either matchings of
$G-\{a,b\}$ and $G-\{c,d\}$, or matchings of $G-\{a,d\}$ and $G-\{b,c\}$, but
never both.  Once again, the cycles and doubled edges are split between the
matchings as described earlier.  Without loss of generality, assume that
paths connect $a$ to $b$ and $c$ to $d$.  As shown earlier, the edge incident
to $a$ must be in the same matching as the edge incident to $b$.  A matching of
$G-\{c,d\}$ may contain both of those edges, but matchings of $G-\{a,d\}, 
G-\{b,c\}$, and $G-\{a,b\}$ cannot.
Likewise, the edges incident to $c$ and $d$ can both belong only to a matching 
of $G-\{a,b\}$.
Thus it is possible for $H$ to be partitioned into matchings of $G-\{a,b\}$
and $G-\{c,d\}$, but not into matchings of $G-\{a,d\}$ and $G-\{b,c\}$.
And just as in the previous paragraph, the partitioning can be done in $2^k$
ways (where $k$ is the number of cycles in $G$).  Thus the number of partitions of $H$ into matchings of $G$
and $G-\{a,b,c,d\}$ is equal to the number of partitions into matchings of
$G-\{a,b\}$ and $G-\{c,d\}$, or of $G-\{a,d\}$ and $G-\{b,c\}$.

Thus we can partition $\Match{G} \times \Match{G-\{a,b,c,d\}}$ and
$(\Match{G-\{a,b\}} \times \Match{G-\{c,d\}}) \cup (\Match{G-\{a,d\}} \times
\Match{G-\{b,c\}})$ into subsets such that the union of each pair of graphs
within the same
subset forms the same graph in $\mathcal{H}$.  Each graph $H \in \mathcal{H}$
corresponds to one subset from each of $\Match{G} \times \Match{G-\{a,b,c,d\}}$
and $(\Match{G-\{a,b\}} \times \Match{G-\{c,d\}})
\cup (\Match{G-\{a,d\}} \times \Match{G-\{b,c\}})$,
and those subsets have equal size.
Thus $\Match{G} \times \Match{G-\{a,b,c,d\}}$
and $(\Match{G-\{a,b\}} \times \Match{G-\{c,d\}})
\cup (\Match{G-\{a,d\}} \times \Match{G-\{b,c\}})$
have the same cardinality, so the relation is proved.
\end{proof}

Before Theorem~\ref{thm:abcd} was known, James Propp proved a special case in
which $a,b,c$, and $d$ form a 4-cycle in $G$; see~\cite{urban}.

\begin{corollary}\label{cor:4-cycle}
Let $a,b,c,d$ be four vertices forming a 4-cycle face in a plane bipartite graph
$G$, joined by edges that we will denote by $ab$, $bc$, $cd$, and $da$.  Then
the proportion $P$ of matchings of $G$ that have an alternating cycle at this
face (i.e., the proportion of matchings of $G$ that either contain edges $ab$
and $cd$ or contain edges $bc$ and $da$) is
\[ P = 2(p(ab)p(cd) + p(bc) p(da)) \]
where $p(uv)$ denotes the proportion of matchings of $G$ that contain the
specified edge $uv$.
\end{corollary}

\begin{proof}
We note that for each edge $uv$ in $G$,
\[ p(uv)=\frac{M(G-\{u,v\})}{M(G)}. \]
The number of matchings of $G$ that contain the alternating cycle at $abcd$
is twice the number of matchings of $G-\{a,b,c,d\}$.  Thus
\[ P = \frac{2 M(G-\{a,b,c,d\})}{M(G)}. \]
Then after multiplying the relation in Theorem~\ref{thm:abcd} by
$2/M(G)^2$, we get our result.
\end{proof}

With this same technique, we can prove similar theorems in which we alter
the membership of $a,b,c,$ and $d$ in $V_1$ and $V_2$.

\begin{theorem}\label{thm:ad-bc}
Let $G=(V_1,V_2,E)$ be a plane bipartite graph in which $|V_1|=|V_2|$.  Let
vertices $a$, $b$, $c$, and $d$ appear in a cyclic order on a face of $G$ (as in
Figure~\ref{fig:ab-cd}, right).  If $a,b \in V_1$ and $c,d \in V_2$, then
\[ M(G-\{a,d\})M(G-\{b,c\}) = M(G)M(G-\{a,b,c,d\}) + M(G-\{a,c\})M(G-\{b,d\}). \]
\end{theorem}

\begin{proof}
The proof of
this relation is similar to that of Theorem~\ref{thm:abcd} with
several differences.  In this case, we show that
$\Match{G-\{a,d\}} \times \Match{G-\{b,c\}}$ and
$(\Match{G} \times \Match{G-\{a,b,c,d\}}) \cup
 (\Match{G-\{a,c\}} \times \Match{G-\{b,d\}})$
have the same cardinality.  The combination of a pair of
matchings from either set produces a graph in the set
$\mathcal{H}$ of graphs on the vertices of $G$ in which all vertices have
degree 2 except for $a,b,c,$ and $d$, which have degree 1.
Now consider a graph $H \in \mathcal{H}$.  If paths
connect $a$ to $b$ and $c$ to $d$, then each path has even length.  The
edges at the ends of each path must go into different matchings.  Thus $H$
can be partitioned into matchings of $G-\{a,d\}$ and $G-\{b,c\}$, or into
matchings of $G-\{a,c\}$ and $G-\{b,d\}$.  Otherwise, if $a$ is connected to
$d$ and $b$ to $c$, then each path has odd length.
Then $H$ can be partitioned into matchings of $G-\{a,d\}$ and $G-\{b,c\}$,
or into matchings of $G$ and $G-\{a,b,c,d\}$.

No matter which ways the path connect, $H$ can always be partitioned into
matchings of $G-\{a,d\}$ and $G-\{b,c\}$.
We can also partition $H$ into either
matchings of $G-\{a,c\}$ and $G-\{b,d\}$, or
matchings of $G$ and $G-\{a,b,c,d\}$, but not both.
Moreover, the number of partitions of $H$
into matchings of $G-\{a,d\}$ and $G-\{b,c\}$
is equal to the number of partitions into matchings
of $G-\{a,c\}$ and $G-\{b,d\}$, or of $G$ and $G-\{a,b,c,d\}$.
Thus $\Match{G-\{a,d\}} \times \Match{G-\{b,c\}}$ and
$(\Match{G} \times \Match{G-\{a,b,c,d\}}) \cup
 (\Match{G-\{a,c\}} \times \Match{G-\{b,d\}})$ have the same cardinality.
\end{proof}

\begin{figure}
\insertfig{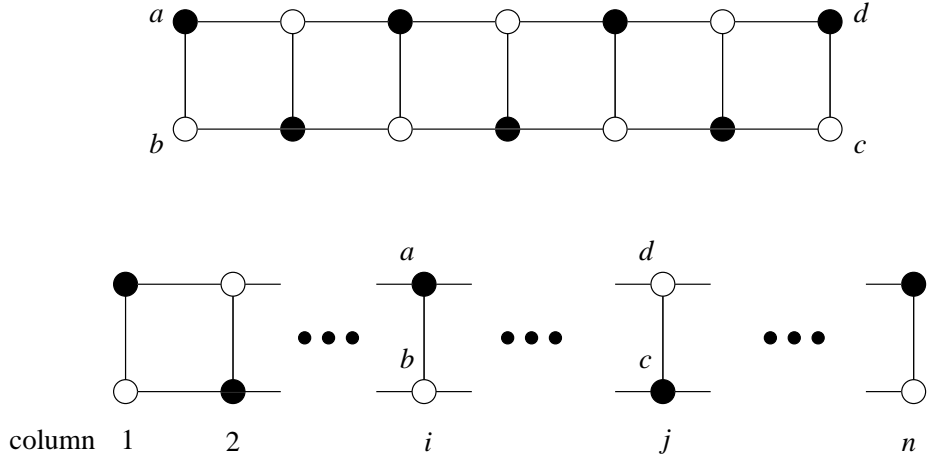}
\caption{Top: $2 \times n$ grid with corners $a,b,c,d.$  Bottom: Same
rectangle, only with $a,b,c,d$ in the middle.}
\label{fig:fibo}
\end{figure}

We show a simple application of Theorems~\ref{thm:abcd} and~\ref{thm:ad-bc} in
which our graphs are $2 \times n$ grids.
From elementary combinatorics, the number of matchings of a $2 \times n$
grid is $F_{n+1}$, where $F_1=F_2=1$, and $F_n=F_{n-1}+F_{n-2}$.  These
theorems lead to a straightforward derivation for some bilinear relations
among the Fibonacci numbers.  Consider a $2 \times n$ rectangle with $a,b,c$
and $d$ being the four corners, as shown in Figure~\ref{fig:fibo}.
Theorems~\ref{thm:abcd} and ~\ref{thm:ad-bc}
produce the relations
\[
F_{n+1}F_{n-1} = \left\{ \begin{array}{ll}
F_n^2 + 1 \cdot 1 & \mbox{for even $n$,} \\
F_n^2 - 1 \cdot 1 & \mbox{for odd $n$,}
\end{array} \right. \]
which could be simplified to $F_{n+1}F_{n-1} = F_n^2 + (-1)^n$.  This is also
known as Cassini's identity.  Another combinatorial proof for this relation 
is found in~\cite{wz86}.

We could go a step further by letting $a,b,c,d$ be somewhere in the middle of
the $2 \times n$ grid.  If $a,b$ are in column $i$, and $c,d$ are in
column $j>i$ (see Figure~\ref{fig:fibo}, bottom), then the relations become
\[ F_{n+1}F_iF_{j-i}F_{n-j+1} =
F_iF_{n-i+1}F_jF_{n-j+1} + (-1)^{j-i-1} (F_iF_{n-j+1})^2, \]
which simplifies to
\[ F_{n+1}F_{j-i} = F_{n-i+1}F_j + (-1)^{j-i-1} F_iF_{n-j+1}. \]

We close this section with two additional relations applicable in situations
in which $V_1$ and $V_2$ have different size.

\begin{theorem}\label{thm:b-acd}
Let $G=(V_1,V_2,E)$ be a plane bipartite graph in which $|V_1|=|V_2|+1$.  Let 
vertices $a$, $b$, $c$, and $d$ appear cyclically on a face of $G$.  If 
$a,b,c \in V_1$ and $d \in V_2$, then
\[ M(G-b)M(G-\{a,c,d\}) = M(G-a)M(G-\{b,c,d\}) + M(G-c)M(G-\{a,b,d\}). \]
\end{theorem}

\begin{theorem}\label{thm:abcd-0}
Let $G=(V_1,V_2,E)$ be a plane bipartite graph in which $|V_1|=|V_2|+2$.  Let 
vertices $a$, $b$, $c$, and $d$ appear cyclically on a face of $G$, and  
$a,b,c,d \in V_1$.  Then
\[ M(G-\{a,c\})M(G-\{b,d\}) = M(G-\{a,b\})M(G-\{c,d\}) + M(G-\{a,d\})M(G-\{b,c\}). \]
\end{theorem}

The proofs for these relations are similar to the proofs for Theorems~\ref{thm:abcd} and~\ref{thm:ad-bc}.

\section{Proof of Aztec Diamond Theorem}\label{sec:adt}

\begin{figure}
\insertfig{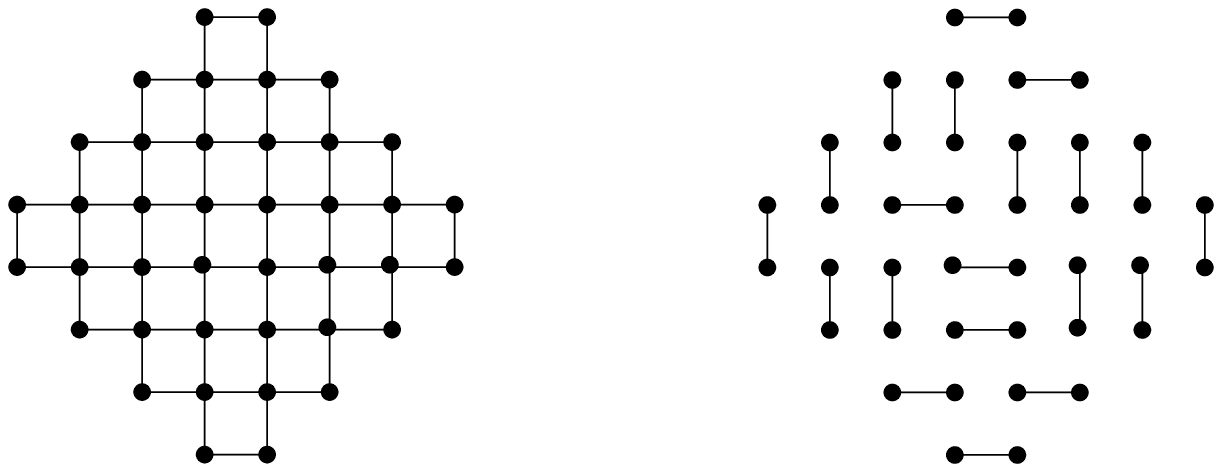}
\caption{Left: Order-4 Aztec diamond graph.  Right: Order-4 Aztec matching.}
\label{fig:aztec-graphs}
\end{figure}

The order-$n$ {\it Aztec diamond graph\/} refers to the graph dual
of the order-$n$ Aztec diamond.  Throughout this proof, an {\it Aztec
matching\/} will mean a matching of an Aztec diamond graph.
Figure~\ref{fig:aztec-graphs} shows the order-4 Aztec diamond graph and an
order-4 Aztec matching.
Thus counting tilings for an Aztec diamond of order $n$ is the same as
counting Aztec matchings of order $n$.

To prove that the number of Aztec matchings of order $n$ is
$2^{n(n+1)/2}$, we need the following recurrence relation.

\begin{proposition}\label{prop:recur-aztec}
Let $T(n)$ represent the number of Aztec matchings of order $n$.  Then
\[ T(n) = \frac{2(T(n-1))^2}{T(n-2)}. \]
\end{proposition}

\begin{proof}
It is sufficient to show that
\[ T(n)T(n-2) = 2(T(n-1))^2. \]

To prove this relation, we show that the number of
ordered pairs $(A,B)$ is twice the number of ordered pairs $(C,D)$, where
$A$, $B$, $C$, and $D$ are Aztec matchings of orders $n$, $n-2$, $n-1$, and
$n-1$, respectively.

We superimpose an Aztec matching $B$ of order $n-2$ with an order-$n$
Aztec matching $A$ so that the matchings are concentric.
Figure~\ref{fig:3-on-5} shows Aztec matchings of orders
3 and 5, and the result of superimposing the two matchings.  In the
combined graph, the white vertices are shared by both the order-3 and
order-5 matchings.  The black vertices are from the order-5 matching only.
Note that some edges are
shared by both matchings.  Note also that each black vertex has degree 1
in the combined graph, whereas each white vertex has degree 2.

\begin{figure}
\insertfig{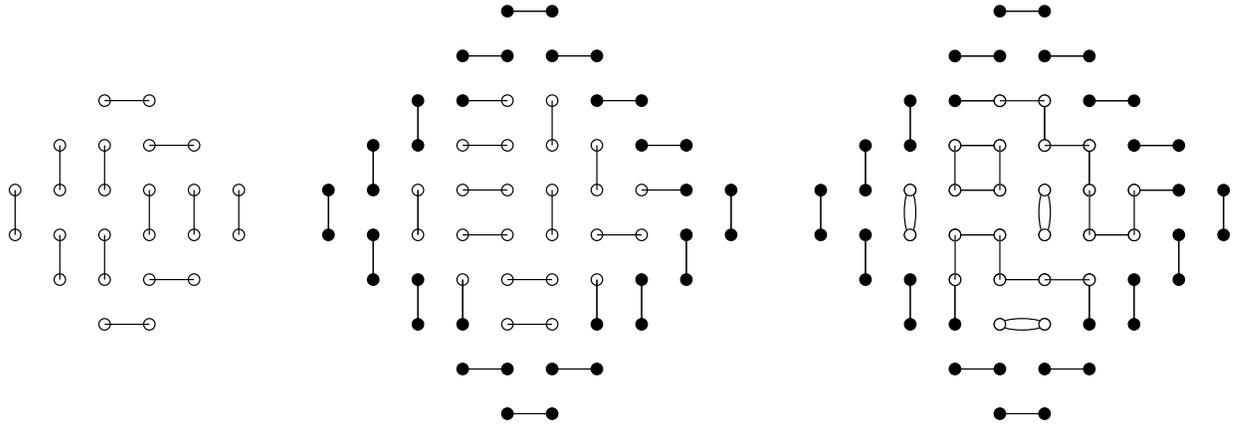}
\caption{Left: Order-3 Aztec matching. Center: Order-5 Aztec matching.
Right: The two matchings combined.}
\label{fig:3-on-5}
\end{figure}

Now consider the two Aztec matchings of order $n-1$ shown in
Figure~\ref{fig:two-ord4}.  Let us call the first and second matchings $C$ and
$D$, respectively.  Figure~\ref{fig:4-on-4} shows the two
possible resulting graphs by superimposing matchings $C$ and $D$ and adding
two extra segments.  The left graph was made
by fitting matching $C$ to the top and matching $D$ to the bottom of the
order-5 diamond, and then adding two side edges.  The graph on the right was
made by fitting matching $C$ to the left and matching $D$ to the right of the
order-5 diamond, and then adding the top and bottom edges.  In both cases,
each of the center vertices has degree 2, and all other vertices have degree 1.
The graphs resemble order-3 Aztec matchings on top of order-5 Aztec matchings.

\begin{figure}
\insertfig{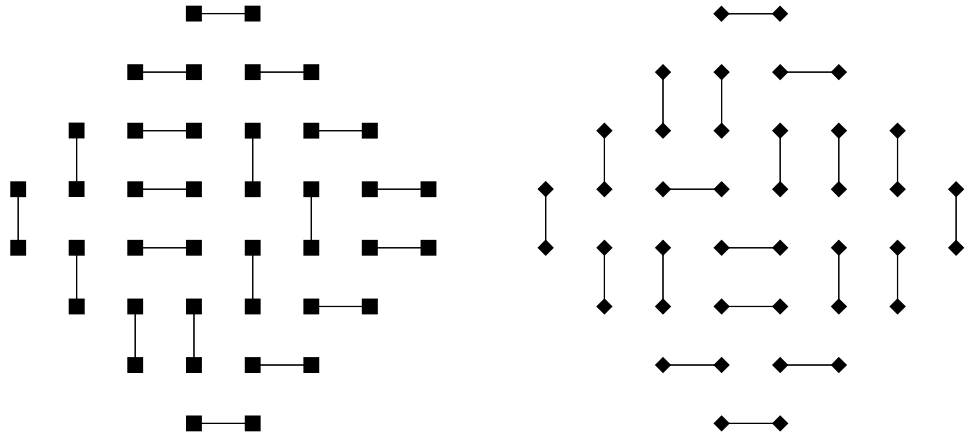}
\caption{Aztec matchings $C$ and $D$ of order 4.}
\label{fig:two-ord4}
\end{figure}

\begin{figure}
\insertfig{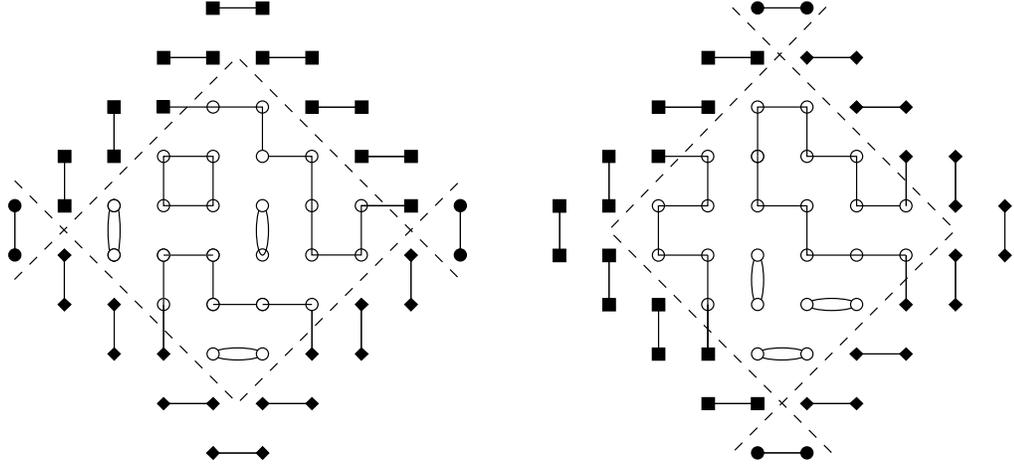}
\caption{The two possible combinations of the matchings $C$ and $D$.}
\label{fig:4-on-4}
\end{figure}

In general, we are given a graph $G$ on the vertices of the order-$n$ Aztec
diamond graph with the following properties:
\begin{enumerate}
\item The inner vertices in $G$ that form an order-$(n-2)$ Aztec diamond
have degree 2.
\item The remaining outer vertices in $G$ have degree 1.
\item The edges of $G$ form cycles, doubled edges, single edges, and
lattice paths of length greater than 1.
\end{enumerate}
Let us call a graph with such properties a {\it doubled Aztec graph.\/}
For each superimposition we have described so far,
the result is a doubled Aztec graph $G$.
We want to show that the number of partitions of $G$ into two Aztec
matchings $A$ and $B$ of orders $n$ and $n-2$ is equal to the number of
partitions of $G$ into two order-$(n-1)$ Aztec matchings $C$ and $D$ (along
with two line segments).  We will show that this number is $2^k$,
where $k$ is the number of
cycles in $G$.  Since $G$ is bipartite, all cycles have even length. These
cycles are contained in the middle common vertices, as they are the only
vertices with degree 2.  Each cycle can then be partitioned so that every
other edge will go to the same subgraph; adjacent edges go to different
subgraphs.  For each cycle, there are two ways to decide which half of the
cycle goes to $A$ or $B$.  Similarly, there are two ways to decide which
half goes to $C$ or $D$. All doubled edges in $G$ are split and shared by
each subgraph.  It remains to show that the other edges must be
partitioned uniquely.

We now label $G$ as shown in Figure~\ref{fig:label-G}.  The vertices whose
degree is 2 are labeled $O$.  The degree-one vertices surrounding the
$O$-vertices are labeled $T$, $U$, $V$, and $W$ such that each side is assigned a
different label. Every vertex on the outer boundary of $G$ is labeled $Y$,
except  for four vertices, one on each corner.  Those four exceptions are
assigned the label ($T, U, V,$ or $W$) of the vertices on the same diagonal.
We have labeled the vertices such that each vertex labeled $Y$ will match
with exactly one vertex
labeled $T, U, V,$ or $W$.  For each label $T, U, V,$ and $W$, exactly one
vertex will not be connected to a $Y$-vertex.  We denote these
special vertices $T'$, $U'$, $V'$, and $W'$.

\begin{figure}
\insertfig{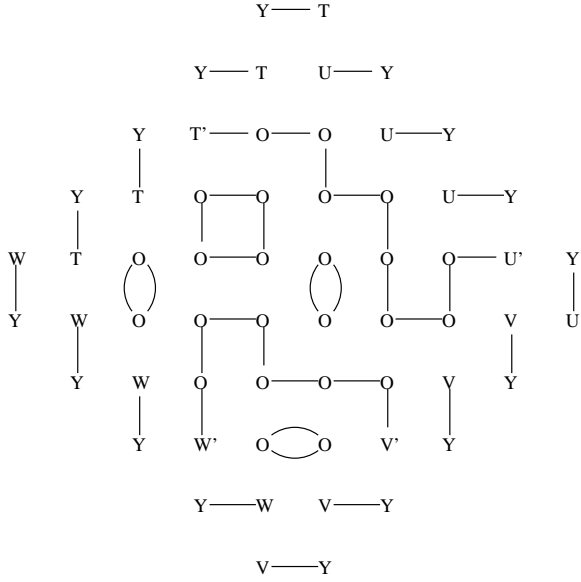}
\caption{Specially labeled graph of doubled Aztec graph $G$.}
\label{fig:label-G}
\end{figure}

In a doubled Aztec graph, there must be paths joining $T'$ to $U'$
and $V'$ to $W'$, or paths joining $T'$ to $W'$ and
$U'$ to $V'$.  However, we cannot have paths going from $T'$ to $V'$
and from $U'$ to $W'$. If such paths existed, then both paths would have to
travel through the $O$-vertices and intersect, thus forcing the degree of some
$O$-vertex to be more than 2.

Let us show that the segments from both ends of a path must belong to the same
subgraph in any partition of $G$.  Let us 2-color the vertices of $G$ black
and white so that black vertices are adjacent to white vertices and vice versa.
The $T$- and $V$- vertices must be the same
color; let us color all the $T$- and $V$-vertices white.  Then the $U$- and
$W$-vertices must be of the other color, which is black.
Therefore, any path from $T$ to $U$, from $U$ to $V$, from $V$ to $W$, or from 
$W$ to $T$ must have odd length since
the path goes from a black to a white vertex.  Thus the segments
from both ends of a path must belong to the same subgraph in any
partition of $G$.

Thus, when we partition
$G$ into matchings $A$ and $B$ of orders $n$ and $n-2$, we must always
place the ending segments into $A$ and determine the rest of the partition
thereafter.  Such a partition always exists.

Next we show that $G$ can be partitioned into two matchings $C$ and $D$
of order $n-1$ along with two additional side edges.
There are two possible ways this partition could be done.
The first is top-bottom: the top diamond contains the $T$- and $U$-vertices,
and the bottom diamond contains the $V$- and $W$-vertices.
The second is left-right: the left diamond contains the $T$- and $W$-vertices,
and the right diamond contains the $U$- and $V$-vertices.

Without loss of generality,
let the paths in $G$ connect $T'$ to $U'$ and $V'$ to
$W'$.  When $G$ is partitioned into two matchings $C$ and $D$, both of order
$n-1$, one matching (say $C$) must have both $T'$ and $U'$,
as they are the ends of the same lattice path.
Thus $C$ is the top Aztec matching containing all $T$- and $U$- vertices
(except for one $U$-vertex on the far right corner).  Vertices $V'$ and $W'$ must
belong to the other matching $D$.  The paths are then partitioned uniquely.
Thus we can partition $G$ into two order-$(n-1)$ Aztec matchings placed
top-bottom (plus two edges on the sides).
However, it is not possible to partition $G$ into two side-by-side Aztec
matchings of order $n-1$ such that one contains the $T$- and $W$-vertices, and
the other contains the $U$- and $V$-vertices.  The reason is that since the left
matching has $W'$, it would then contain $V'$. The latter cannot
happen, since $V'$ is in the other matching.

Hence each doubled Aztec graph can be partitioned into two
order-$(n-1)$ Aztec matchings in one way (top-bottom) or the other
(left-right), but never both.  The partition of the paths is uniquely
determined.

The number of ways to combine Aztec matchings of orders $n$ and $n-2$ is
$T(n)T(n-2)$, while the number of ways to combine two order-$(n-1)$ matchings
is $2T(n-1)^2$.  Each combination becomes a doubled Aztec graph, so the
relation is proved.
\end{proof}

There are 2 ways to tile an order-1 Aztec diamond, and 8 ways to tile
an order-2 Aztec diamond.  Having proved the recurrence relation, we can
now compute the number of tilings of an Aztec diamond of order $n$.  The
following result is easily proved by induction on $n$:

\begin{theorem}[Aztec Diamond Theorem]\label{thm:aztec-diamond}
The number of tilings of the order-$n$ Aztec diamond is $2^{n(n+1)/2}$.
\end{theorem}

\section{Regions with holes}\label{sec:holes}

\subsection{Placement Probabilities}

We can use graphical condensation to derive recurrence relations for placement
probabilities of dominoes in tilings of Aztec diamonds.  Let domino $D$
be a specified pair of adjacent squares in an Aztec diamond.
The {\it placement
probability\/} of $D$ in an order-$n$ Aztec diamond is the
probability that $D$ will appear in a tiling of the order-$n$ Aztec diamond,
given that all tilings are equally likely.

Placement probabilities are of interest in the study of random tilings.
If we look at a random tiling of an Aztec diamond of large order, we notice
four regions in which the dominoes form a brickwork pattern, and a central
circular region where dominoes are mixed up.
The placement probability of any domino
at the center of the diamond will be near 1/4.  However, in the top corner,
dominoes which conform to the brickwork will have probabilities near 1.
All other dominoes in this corner would have probabilities near 0.
For proofs of these assertions, see~\cite{cep}.

We could calculate the placement probability of a domino with the following
steps.  First, we replace the domino with a two-square hole in the
Aztec diamond.
Then we compute the number of tilings of that diamond with the hole.
Finally, we divide it by the number of tilings of the (complete) Aztec
diamond.

\begin{figure}
\insertfig{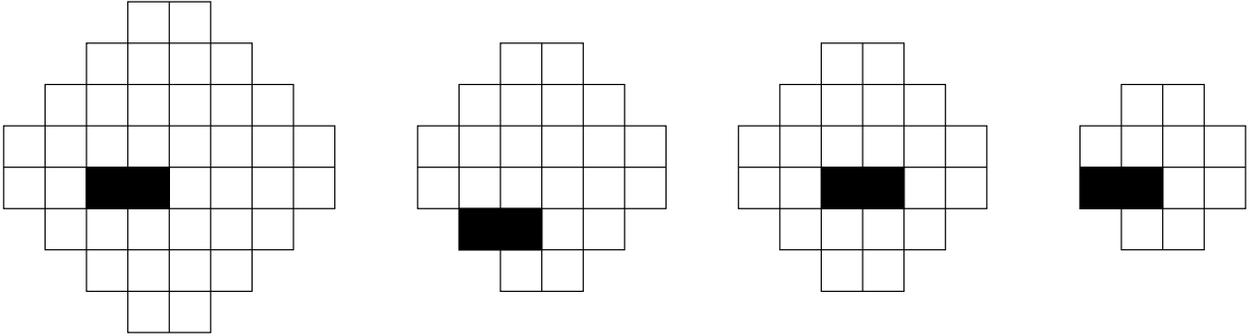}
\caption{Regions from left: $A_n-D$, $A_{n-1}-D_{\rm down}$, $A_{n-1}-D_{\rm right}$,
$A_{n-2}-D$.}
\label{fig:holey-diamonds}
\end{figure}

We can express the number of tilings of the order-$n$ Aztec diamond with
the hole at $D$ in terms of tilings of lower-order Aztec diamonds with
holes.  But first, let us introduce some notation.
We will let $A_{n}-D$ stand for the
order-$n$ Aztec diamond with domino $D$ missing.  The dominoes $D_{\rm up}$,
$D_{\rm down}$, $D_{\rm left}$, and $D_{\rm right}$ will represent dominoes shifted up,
down, left, and right by a square relative to $D$ in the Aztec diamond.
Then $A_{n-1}-D_{\rm up}$ is the order-$(n-1)$ diamond such that when it is placed
concentrically with $A_n-D$, the hole of $A_{n-1}-D_{\rm up}$ will match up with
$D_{\rm up}$.  The regions $A_{n-1}-D_{\rm down}$ and so forth represent similar Aztec
diamonds with domino holes.  Finally, $A_{n-2}-D$ is the order-$(n-2)$
Aztec diamond such that when $A_{n-2}-D$ is placed directly over $A_n-D$, domino
$D$ is missing.  See Figure~\ref{fig:holey-diamonds} for examples.  (In case
$D_{\rm up}$, etc. lies outside $A_{n-1}$, the region $A_{n-1}-D_{\rm up}$ will
not be defined.)

By using graphical condensation, we can relate the number of tilings of these
Aztec diamonds with holes:
\[T(A_n-D) T(A_{n-2}-D) = T(A_{n-1}-D_{\rm up}) T(A_{n-1}-D_{\rm down})
                        + T(A_{n-1}-D_{\rm left}) T(A_{n-1}-D_{\rm right}).\]

We also have the following relation, which relates
numbers of tilings of Aztec diamonds:
\[T(A_n) T(A_{n-2}) = 2 T(A_{n-1}) T(A_{n-1}).\]

We can then derive a relation among placement probabilities of dominoes in
Aztec diamonds of orders $n$, $n-1$, and $n-2$.  When we divide the first
relation by the second, we get
\[P(A_n, D) P(A_{n-2}, D) = \frac{P(A_{n-1}, D_{\rm up}) P(A_{n-1},D_{\rm down})
+ P(A_{n-1},D_{\rm left}) P(A_{n-1},D_{\rm right})}{2},\]
where $P(R, D)$ is the placement probability on domino $D$ in region $R$.
The probability $P(R,D)$ was computed by dividing $T(R-D)$ by $T(R)$.

\subsection{Holey Aztec Rectangles}

Another application of graphical condensation deals with regions called
holey Aztec rectangles.  A {\it holey Aztec rectangle} is a region similar
to an Aztec diamond, except that the boundary of an $n$-by-$(n+1)$ holey
Aztec rectangle consists of diagonals of length $n$, $n+1$, $n$, and $n+1$.
In addition, to maintain the balance of squares of different parity so that
the region can be tiled, a square is removed from its interior.  Problems 9
and 10 in~\cite{enum} ask to enumerate tilings of a holey Aztec rectangle
with a square removed in the center or adjacent to the center square, depending
on the parity of $n$.

\begin{figure}
\insertfig{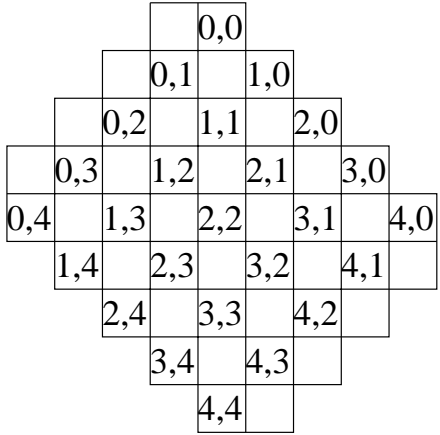}
\caption{Labeling of squares in a holey Aztec rectangle.  By removing square
$(a,b)$ we obtain $R_{4,a,b}$.}
\label{fig:labeled-har}
\end{figure}

Let us label some of the squares in an Aztec rectangle as shown in
Figure~\ref{fig:labeled-har}.  We label a square only if
the region becomes tileable after deleting that square.
We let $R_{n,a,b}$ represent the $n$-by-$(n-1)$ Aztec
rectangle whose square $(a,b)$ has been deleted.   Then we can apply our
technique and come up with a theorem which relates the numbers of tilings
among holey Aztec Rectangles.

\begin{theorem}\label{thm:har}
Let $T(R)$ stand for the number of tilings of a region $R$. Then for $a$, $b$
between $1$ and $n-1$, the number of tilings of $R_{n,a,b}$ is expressed in the
following relation:
\[T(R_{n,a,b}) =
\frac{T(R_{n-1,a,b-1})T(R_{n-1,a-1,b})
+T(R_{n-1,a,b})T(R_{n-1,a-1,b-1})}{T(R_{n-2,a-1,b-1})}.\]
\end{theorem}

\begin{proof}
The proof is very similar to that of Proposition~\ref{prop:recur-aztec}.
Instead of superimposing an order-$n$ Aztec matching
on an order-$(n-2)$ Aztec matching,
we superimpose $R_{n-2,a-1,b-1}$ on top of $R_{n,a,b}$ so that the holes
align to the same spot.  Given the graph $G$ resulting from the superimposition,
we can partition it into two $(n-1)$-by-$n$ holey Aztec rectangles.  The
partition can be done either left-right or top-bottom, but only one or the
other.  The left-right rectangles are isomorphic to $R_{n-1,a-1,b}$ and
$R_{n-1,a,b-1}$.  The top-bottom rectangles are isomorphic to
$R_{n-1,a,b}$ and $R_{n-1,a-1,b-1}$.
\end{proof}

Another relation can be proven for the case in which the hole
is on the edge of the rectangle:

\begin{theorem}\label{thm:har2}
If $1 \leq a \leq n$, then
\[ T(R_{n,a,0}) =
\frac{T(R_{n-1,a,0})T(A_n)+T(R_{n-1,a-1,0})T(A_n)}{T(A_{n-1})} \]
where $A_n$ is the Aztec diamond of order $n$.
\end{theorem}

\begin{figure}
\insertfig{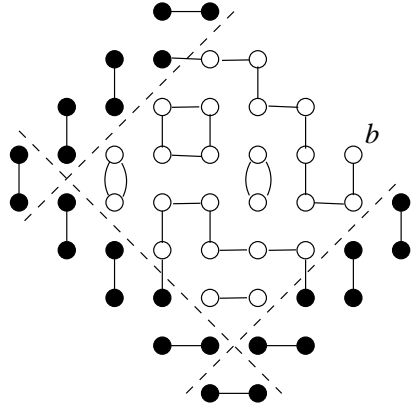}
\caption{Combined matchings of an order-3 Aztec diamond graph and $R_{4,3,0}$.}
\label{fig:4-5-double}
\end{figure}

\begin{proof}
As an example, Figure~\ref{fig:4-5-double} shows a matching of an order-3
Aztec diamond graph
(which is shown in white vertices) on a matching of $R_{4,3,0}$ (which is
missing the vertex $b$).  The relation is derived in a manner analogous to
Theorem~\ref{thm:b-acd}.
\end{proof}

\subsection{``Pythagorean'' regions}

\begin{figure}
\insertfig{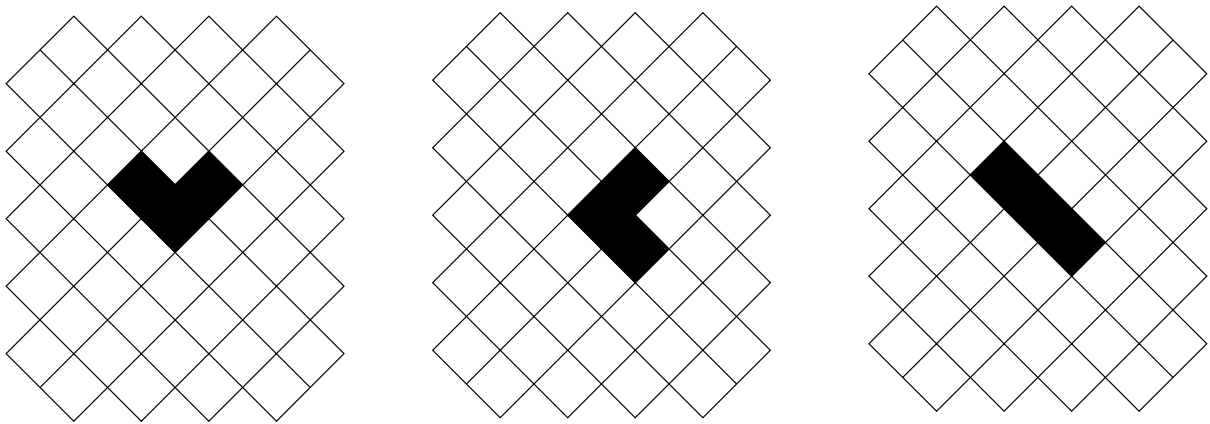}
\caption{From left to right: $R_4-t_1$, $R_4-t_2$, and $R_4-t_3$.}
\label{fig:pythagorean}
\end{figure}

We can derive one more relation as a corollary to Theorem~\ref{thm:abcd-0}.
Let $R_n$ be an $n \times (n+1)$ Aztec rectangle, where $n$ is
even.  Let $t_1, t_2$, and $t_3$ be (overlapping) trominoes in $R_n$, each of
which contain the center square and two squares adjacent to it.  Trominoes
$t_1$ and $t_2$ are L-shaped, while $t_3$ is straight.  Let $t_1$ point to a
side of length $n$, and $t_2$ point to a side of length $n+1$.  (See 
Figure~\ref{fig:pythagorean}.)  Then
\[ T(R_n-t_1)^2 + T(R_n-t_2)^2 = T(R_n-t_3)^2. \]
In other words, we have a Pythagorean relation among the number of tilings of
these regions!  The proof of this relation is to set $G$ to be $R_n$ minus
the center square, let $a,b,c,d$ be squares adjacent to the center, and then
apply Theorem~\ref{thm:abcd-0}.

\begin{figure}
\insertfig{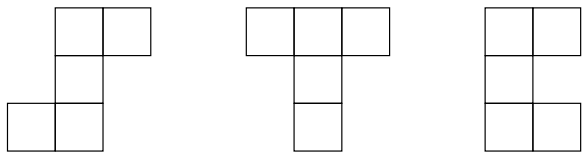}
\caption{Pentominoes missing from rectangular regions}
\label{fig:pentominoes}
\end{figure}

The reader may also like to puzzle over a similar ``Pythagorean'' relation
among the numbers of tilings of
rectangular ({\it not\/} Aztec rectangular) regions
in which each region has a pentomino hole in its center.
The pentominoes are shown in Figure~\ref{fig:pentominoes}.

\section{Weighted matchings of Planar Bipartite Graphs and Aztec Diamonds}
\label{sec:weighted}

\subsection{Weighted Planar Bipartite Graphs}

We can generalize the enumerative relations proved in section~\ref{sec:pbg}
to cover weighted planar bipartite graphs.
Given a graph $G$, we can assign a weight to each edge to form a {\it
weighted graph\/}.  The {\it weight\/} of any subgraph $H$ of $G$
is the product of the weights of all the edges in $G$
(in the case where $H$ is a multigraph,
each edge-weight contributes with exponent equal to
the multiplicity of the associated edge in $H$);
e.g., the weight of a matching of $G$
is the product of the weights of each edge in that matching.
We denote the weight of $G$ itself by $w(G)$.
We also define the {\it weighted sum\/} $W(G)$ of $G$ to be
the sum of the weights of all possible matchings on $G$.

We can now state and prove a weighted version of Theorem~\ref{thm:abcd}:

\begin{theorem}\label{thm:wtd-abcd}
Let $G=(V_1,V_2,E)$ be a weighted plane bipartite graph in which
$|V_1|=|V_2|$.  Let vertices $a$, $b$, $c$, and $d$ appear on a face of $G$,
in that order.  If $a, c \in V_1$ and $b,d \in V_2$, then
\[ W(G) W(G-\{a,b,c,d\}) = W(G-\{a,b\})W(G-\{c,d\}) + W(G-\{a,d\})W(G-\{b,c\}). \]
\end{theorem}

\begin{proof}
The proof essentially follows that of Theorem~\ref{thm:abcd}, except that we
must now account for the weights.  Let $\mathcal{H}$ be the set of
graphs on the vertices of $G$ in which vertices $a$, $b$, $c$, and $d$
have degree 1, all other vertices have degree 2, and doubled edges are
permitted.  Let $H$ be a graph in $\mathcal{H}$.  As before, $H$ may be
partitioned into two matchings $M_1$ and $M_2$ with these possibilities:
\begin{enumerate}
\item $(M_1,M_2) \in \Match{G} \times \Match{G-\{a,b,c,d\}}$.
\item $(M_1,M_2) \in \Match{G-\{a,b\}} \times \Match{G-\{c,d\}}$.
\item $(M_1,M_2) \in \Match{G-\{a,d\}} \times \Match{G-\{b,c\}}$.
\end{enumerate}

As we have seen before, $H$ can always be partitioned in choice 1, and also in
either choice 2 or choice 3 (but not both).  The number of possible partitions
is $2^k$, where $k$ is the number of cycles in $H$.  So
\begin{eqnarray*}
W(G)W(G-\{a,b,c,d\}) &=& \sum_{H \in \mathcal{H}} 2^{k(H)} w(H) \\
                &=& W(G-\{a,b\})W(G-\{c,d\}) + W(G-\{a,d\})W(G-\{b,c\}),
\end{eqnarray*}
where $k(H)$ is the number of cycles in $H$.
\end{proof}

Similar relations can be generalized from
Theorems~\ref{thm:ad-bc},~\ref{thm:b-acd}, and~\ref{thm:abcd-0}:

\begin{theorem}\label{thm:wtd-ad-bc}
Let $G=(V_1,V_2,E)$ be a weighted plane bipartite graph
in which $|V_1|=|V_2|$.  Let
vertices $a$, $b$, $c$, and $d$ appear on a face of $G$, in that order (as in
Figure~\ref{fig:ab-cd}, right).  If $a,b \in V_1$ and $c,d \in V_2$, then
\[ W(G-\{a,d\})W(G-\{b,c\}) =
W(G)W(G-\{a,b,c,d\}) + W(G-\{a,c\})W(G-\{b,d\}). \]
\end{theorem}

\begin{theorem}\label{thm:wtd-b-acd}
Let $G=(V_1,V_2,E)$ be a weighted plane bipartite graph in which
$|V_1|=|V_2|+1$.  Let vertices $a$, $b$, $c$, and $d$ appear on a face of $G$,
in that order.  If $a,b,c \in V_1$ and $d \in V_2$, then
\[ W(G-b)W(G-\{a,c,d\}) = W(G-a)W(G-\{b,c,d\}) + W(G-c)W(G-\{a,b,d\}). \]
\end{theorem}

\begin{theorem}\label{thm:wtd-abcd-0}
Let $G=(V_1,V_2,E)$ be a weighted plane bipartite graph in which
$|V_1|=|V_2|+2$.  Let vertices $a$, $b$, $c$, and $d$ appear on a face of $G$,
in that order,and  $a,b,c,d \in V_1$.  Then
\[ W(G-\{a,c\})W(G-\{b,d\}) =
W(G-\{a,b\})W(G-\{c,d\}) + W(G-\{a,d\})W(G-\{b,c\}). \]
\end{theorem}

\subsection{Weighted Aztec Diamonds}

\begin{figure}
\insertfig{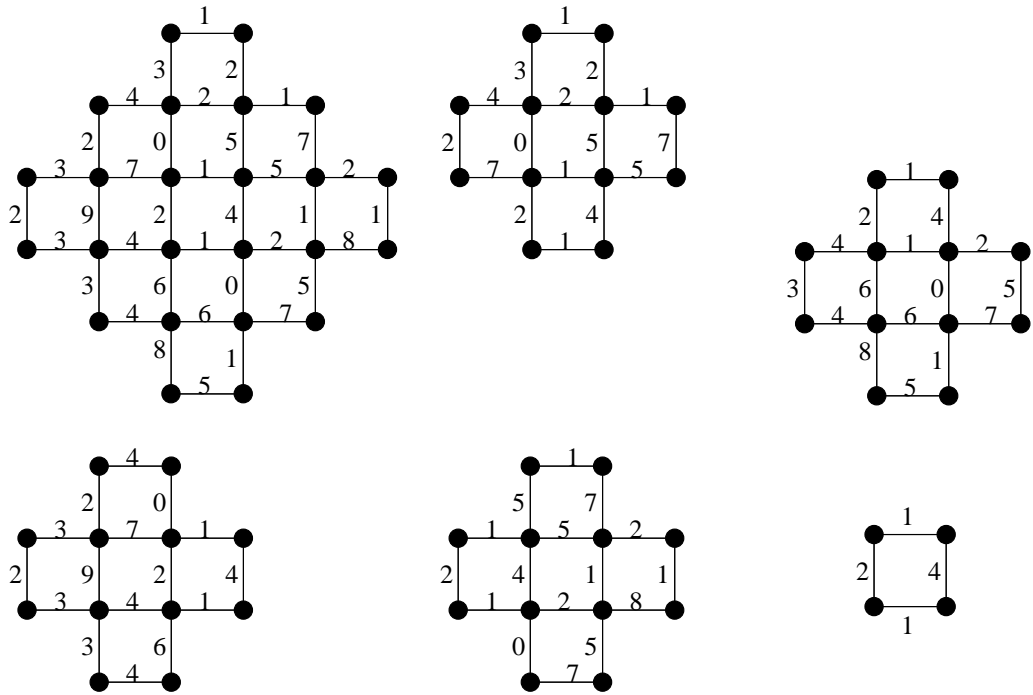}
\caption{Top Row: A weighted Aztec diamond graph $A$, followed by $A_{\rm top}$
and $A_{\rm bottom}$. Bottom Row: $A_{\rm left}$, $A_{\rm right}$, and $A_{\rm middle}$.}
\label{fig:weighted}
\end{figure}

Consider a weighted Aztec diamond graph $A$ of order $n$.
Define $A_{\rm top}$ to be
the upper order-$(n-1)$ Aztec sub-diamond along with its corresponding edge
weights in $A$.  Similarly, we can refer to the bottom, left, and
right subgraphs of $A$ as $A_{\rm bottom}$, $A_{\rm left}$, and $A_{\rm right}$, which are
all order $n-1$ Aztec sub-diamonds.  Finally, let $A_{\rm middle}$ be the
inner order-$(n-2)$ Aztec diamond within $A$.  Figure~\ref{fig:weighted}
shows an Aztec diamond and its five sub-diamond graphs.

It turns out that the superimposition technique can also be used to establish
an identity for weighted Aztec diamond graphs.  The following theorem shows
how the weighted sum of a weighted Aztec
diamond can be expressed in terms of the weighted sums of the
subdiamonds and a few edge weights.

\begin{theorem}\label{thm:wtd-aztec}
Let $A$ be a weighted Aztec diamond of order $n$.
Also let $t$, $b$, $l$, and $r$ be the weights of the top, bottom, left, and
right edges of $A$, respectively.  Then
        \[ W(A) = \frac{l\cdot r\cdot W(A_{\rm top})\cdot W(A_{\rm bottom})
+ t\cdot b\cdot W(A_{\rm left})\cdot W(A_{\rm right})}{W(A_{\rm middle})}. \]
\end{theorem}

\begin{proof}
This proof is very similar to Proposition~\ref{prop:recur-aztec},
except that we must fill in the details concerning the weights.
Indeed, we want to show that
\begin{equation}\label{eq:wtd-az}
W(A)\cdot W(A_{\rm middle})
= l\cdot r\cdot W(A_{\rm top})\cdot W(A_{\rm bottom})
+ t\cdot b\cdot W(A_{\rm left})\cdot W(A_{\rm right}).
\end{equation}

We have seen how a doubled Aztec graph $G$ of order $n$ can be decomposed into
subgraphs in two of three following ways:
\begin{enumerate}
\item (Big-small) Two Aztec matchings of orders $n$ and $n-2$.
\item (Top-bottom) Top and bottom Aztec matchings of order $n-1$, plus the
left and right edges.
\item (Left-right) Left and right Aztec matchings of order $n-1$, plus the
top and bottom edges.
\end{enumerate}
As we know, $G$ can always be decomposed via Big-small and
by either Top-bottom or Left-right (but not both).  The number of possible
decompositions by either method is $2^k$, where $k$ is the number of cycles
in $G$.

Each edge in $G$ becomes a part of exactly one of the subgraphs.  Therefore,
the product of the weights of the subgraphs will always equal to the weight of
$G$, since each edge weight is multiplied once.

Recall that $W(A)$ is the sum of the weights of all possible matchings
on $A$.  Then
        \[ W(A)\cdot W(A_{\rm middle}) = \sum_{G} 2^{k(G)}w(G) \]
where $G$ ranges over all doubled Aztec graphs of order $n$, and $k(G)$ is
the number of cycles in $G$.  Each term in the sum represents the weight
of $G$ multiplied by the number of ways to partition $G$ via Big-small.
Each partition is accounted for in $W(A)W(A_{\rm middle}).$  Similarly, we
also have
\[l\cdot r\cdot W(A_{\rm top})\cdot W(A_{\rm bottom})
+ t\cdot b\cdot W(A_{\rm left})\cdot W(A_{\rm right}) = \sum_{G} 2^{k(G)}w(G). \]

Thus both sides of Equation~\ref{eq:wtd-az} are equal to a common
third quantity, so the relation is proved.
\end{proof}

\begin{figure}
\insertfig{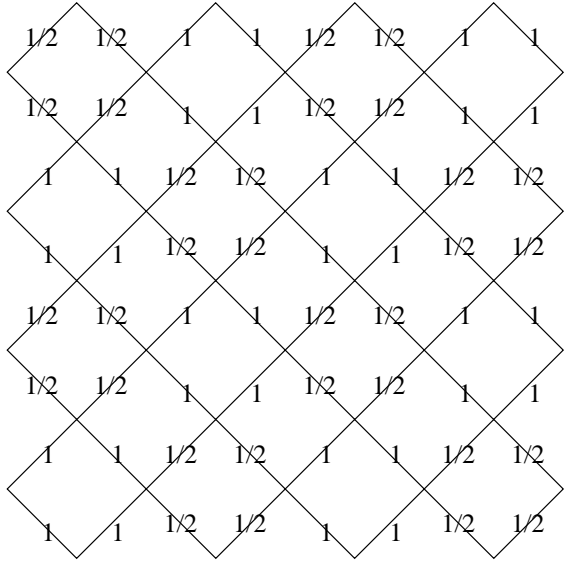}
\caption{A fortress-weighted Aztec diamond, rotated by 45 degrees.}
\label{fig:fort}
\end{figure}

Theorem~\ref{thm:wtd-aztec} may be used to find the weighted
sum of a {\it fortress-weighted Aztec diamond\/}.  Imagine rotating an
Aztec diamond graph by 45 degrees and then partitioning the edges of the
graph into {\it cells\/}, or sets of four edges forming a cycle.  In a
fortress-weighted Aztec diamond, there are two types of cells: (1) cells
whose edges are weight 1, and (2) cells whose edges are weight 1/2.
Cells with edge-weights of 1/2 are adjacent to cells with edge-weights of 1.
(See Figure~\ref{fig:fort}).

There are three kinds of fortress-weighted Aztec diamonds:
\begin{enumerate}
\item The order $n$ is odd, and all edges in the corner cells have weight 1.
\item The order $n$ is odd, and all edges in the corner cells have weight 1/2.
\item The order $n$ is even, and two opposite corners have edges weighted 1/2,
and the other two corners have edges weighted 1.
\end{enumerate}

Let $A_n$, $B_n$, and $C_n$ stand for the weighted sums of the these diamonds,
respectively.  We then use Theorem~\ref{thm:wtd-aztec} to establish relations
among $A_n$, $B_n$, and $C_n$.  They are

\begin{eqnarray*}
A_{2k+1} &=& \frac{1\cdot 1\cdot C_{2k}\cdot C_{2k}
  + 1\cdot 1\cdot C_{2k}\cdot C_{2k}}{A_{2k-1}}
  = \frac{2{C_{2k}}^2}{A_{2k-1}}, \\
B_{2k+1} &=& \frac{\frac{1}{2}\cdot \frac{1}{2}\cdot C_{2k}\cdot C_{2k}
  + \frac{1}{2}\cdot \frac{1}{2}\cdot C_{2k}\cdot C_{2k}}{B_{2k-1}}
  = \frac{\frac{1}{2} {C_{2k}}^2}{B_{2k-1}}, \\
C_{2k} &=& \frac{\frac{1}{2}\cdot \frac{1}{2}\cdot A_{2k-1}\cdot A_{2k-1}
  + 1\cdot 1\cdot B_{2k-1}\cdot B_{2k-1}}{C_{2k-2}}
  = \frac{\frac{1}{4} {A_{2k-1}}^2 + {B_{2k-1}}^2}{C_{2k-2}}.
\end{eqnarray*}

{}From these relations, we can easily prove by induction that for odd $k$,
\begin{eqnarray*}
A_{2k+1} &=& (5/4)^{k(k+1)}, \\
B_{2k+1} &=& (5/4)^{k(k+1)}, \\
C_{2k} &=& (5/4)^{k^2}.
\end{eqnarray*}
For even $k$,
\begin{eqnarray*}
A_{2k+1} &=& 2(5/4)^{k(k+1)}, \\
B_{2k+1} &=& \frac{(5/4)^{k(k+1)}}{2}, \\
C_{2k} &=& (5/4)^{k^2}.
\end{eqnarray*}

\begin{figure}
\insertfig{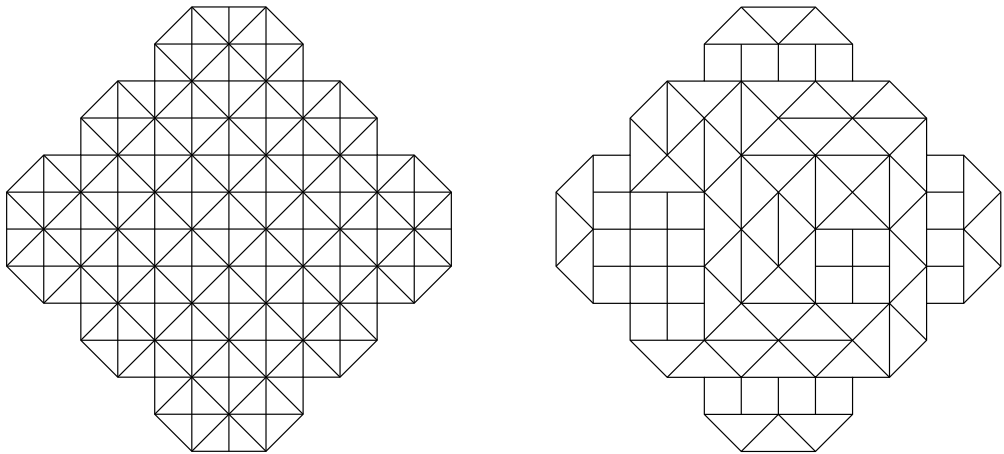}
\caption{A fortress and sample diabolo tiling.}
\label{fig:fortress}
\end{figure}

The importance of the fortress-weighted Aztec diamond comes from the problem
of computing the number of diabolo tilings for a fortress.
A {\it diabolo\/} is either an isosceles right triangle or a square, formed by
joining
two smaller isosceles right triangles.  A {\it fortress\/} is a diamond shaped
region that is made up of isosceles right triangles and can be tiled by
diabolos.  A fortress and a sample tiling by
diabolos are shown in Figure~\ref{fig:fortress}.

To transform a fortress graph into a weighted Aztec diamond graph, we must use
a method called {\it urban renewal.\/} This technique is explained
in~ \cite{urban} along with proofs and applications.  
In~\cite{urban}, the transformation is described for the fortress, and the
number of tilings for the fortress would be the weighted sum of the
fortress-weighted diamond times some power of 2. Thus, graphical condensation,
in combination with this known result about enumeration of fortresses,
provides a very simple way to derive the formulas for the number of fortress 
tilings, first proven by Bo-Yin Yang~\cite{yang}.

A different sort of weighting scheme allows us to apply graphical
condensation to count domino tilings of ordinary (non-Aztec!) rectangles.
Every rectangle $R$ of even area can be imbedded in some Aztec diamond
$A$ of order $n$ (with $n$ sufficiently large) in such a fashion that
the complement $A \setminus R$ (the portion of $A$ that is not covered
by $R$) can be tiled by dominoes $d_1,\dots,d_M$.
For any such tiling of $A \setminus R$, we can define a weighting of
the Aztec diamond graph of order $n$ with the property that each matching
of $A$ has weight 1 if the associated tiling of $A$ contains all the
dominoes $d_1,\dots,d_M$ and weight 0 otherwise.  (Specifically, assign
weight 1 to every edge that corresponds to one of the dominoes $d_1,\dots,d_M$
or to a domino that lies entirely inside R, and weight 0 to every other
edge.)  Then the sum of the weights of the matchings of the weighted
Aztec diamond graph equals the number of tilings of $R$.

\section{Plane Partitions}\label{sec:pp}

A plane partition is a finite array of integers such that each row and column
is a weakly decreasing sequence of nonnegative integers.  If we represent each
integer $n$ in the plane partition as a stack of $n$ cubes, then the plane
partition is a collection of cubes pushed into the corner of a box.  When this
collection of cubes is viewed at a certain angle, these cubes will appear as a
rhombus tiling of a hexagon.

In 1912, Percy MacMahon~\cite{M12} published a proof of a generating function
that enumerates plane partitions that fit in a box $\bx{r,s,t}$ with dimensions
$r \times s \times t$.

\begin{theorem}\label{thm:pp}
Define $P(r,s,t)$ as the generating function for plane partitions that fit in
$\bx{r,s,t}$.  Then
\[ P(r,s,t)=\prod_{i=1}^r \prod_{j=1}^s \frac{1-q^{i+j+t-1}}{1-q^{i+j-1}}. \]
\end{theorem}

Other proofs of this theorem have been published by Carlitz~\cite{C67},
and Gessel and Viennot~\cite{GV85}.

In this section, we will prove MacMahon's formula with the help of graphical
condensation.  Using graphical condensation, we derive a relation that
enables us to prove MacMahon's formula by induction on $r+s+t$.

\begin{theorem}\label{thm:wtd-hex-1}
\[P(r+1,s+1,t)P(r,s,t)
= q^t P(r,s+1,t) P(r+1,s,t) + P(r+1,s+1,t-1) P(r,s,t+1). \]
\end{theorem}

\begin{proof}
Let us take the dual graph of a hexagonal region of triangles in
which $r$ is the length of the bottom right side,
$s$ is the length of the bottom left side, and $t$ is the height of the
vertical side.  In this dual graph,
all edges that are not horizontal are weighted 1.  The horizontal edges are
weighted as follows: the $r$ edges along the bottom right diagonal are each
weighted 1.  On the next diagonal higher up, each edge is weighted $q$, and
the weights of the edges on each subsequent diagonal are $q$ times the weights
of the previous diagonal.  Thus the range of weights should be from $1$ to
$q^{s+t-1}$.  (See Figure~\ref{fig:wtd-hex}.)  Call this weighted graph
$H(r,s,t)$.

\begin{figure}
\insertfig{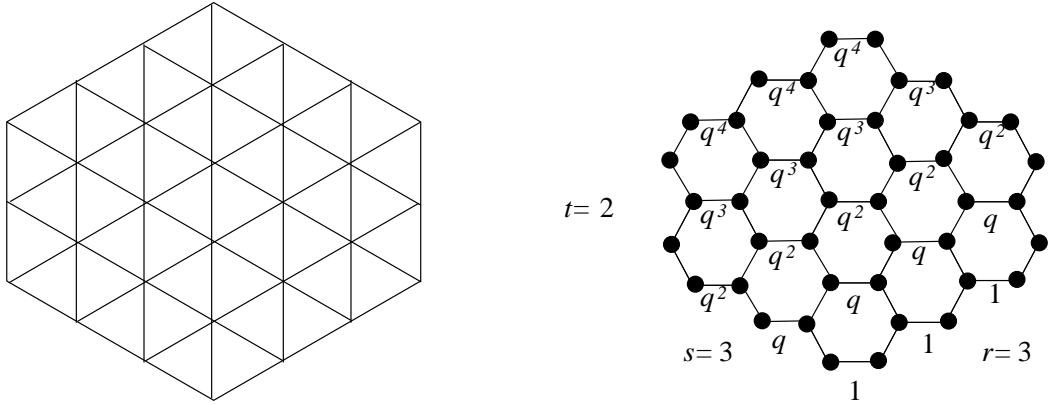}
\caption{Weighting scheme for $H(3,3,2)$.  All unlabeled edges have weight 1.}
\label{fig:wtd-hex}
\end{figure}

This weighting scheme is specifically designed so that, if a matching consists
of the bottom edge (weighted $q^i$) and two other edges of a 6-cycle, then by
replacing those edges with the other three edges, we have dropped the 
$q^i$-weighted edge in favor of the $q^{i+1}$-weighted edge.  (See
Figure~\ref{fig:wtd-6-cycle}.)  The matching would then gain a factor of $q$,
resembling the action of
adding a new block (weighted $q$) to a plane partition.  The minimum weight of
a matching of this graph is $q^{rs(s-1)/2}$, corresponding to the $rs$
horizontal edges that would make up the ``floor'' of the empty plane partition.
The weighted sum of the graph is therefore $q^{rs(s-1)/2}P(r,s,t)$.

\begin{figure}
\insertfig{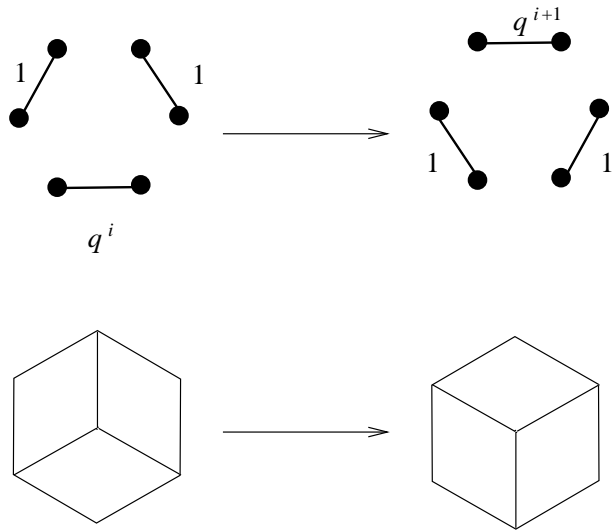}
\caption{A transition representing the addition of another cube to a Young
diagram.}
\label{fig:wtd-6-cycle}
\end{figure}

Now the proof of this relation is very similar to the proofs of
Theorem~\ref{thm:abcd} and Proposition~\ref{prop:recur-aztec}.
We superimpose the two weighted hexagonal graphs
$H(r+1,s+1,t)$ and $H(r,s,t)$ such that the bottom edge common to sides $r$
and $s$ of $H(r,s,t)$ coincides with the bottom edge of $H_{r-1,s-1,t}$.  The
two hexagons completely overlap except for four outer strips of triangles from
$H(r+1,s+1,t)$.  Let us number these strips 1, 2, 3, and 4.  (See
Figure~\ref{fig:wtd-hex-1}.)
When we superimpose the two matchings in the manner
described above, we get once again a collection of cycles, doubled edges,
single edges, and two paths.  Each vertex inside $H(r,s,t)$ has degree 2, and
each vertex in the four outer strips has degree 1.

\begin{figure}
\insertfig{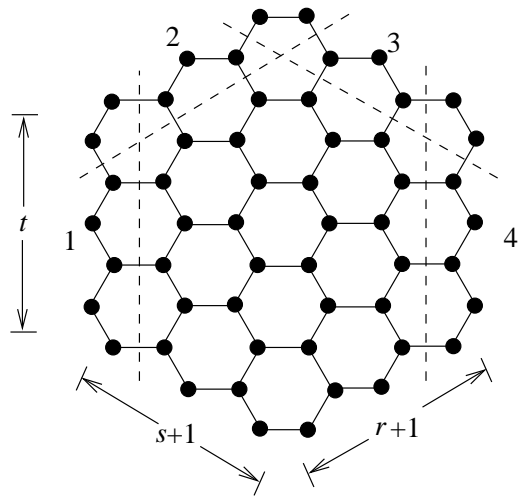}
\caption{Graphs $H(r+1,s+1,t)$ and $H(r,s,t)$ overlapping.  The outer strips
 are numbered 1,2,3,4.}
\label{fig:wtd-hex-1}
\end{figure}

Within each strip, all but one of the vertices are matched with each other.
Those four unmatched vertices are
the endpoints of the two paths.  If one path runs between vertices on strips
1 and 2, and the other runs between vertices on strips 3 and 4, then the
collection can be partitioned into matchings of the duals of $H(r,s+1,t)$
and $H(r+1,s,t)$, plus the edge on the corner of strips 2 and 3 (of weight
$q^{s+t-1}$).  The graph $H(r,s+1,t)$ lacks
strips 3 and 4, while $H(r+1,s,t)$ is the graph without
strips 1 and 2.  Alternatively, if the paths run from strip 1
to strip 4, and from strip 2 to strip 3, then the collection
can be partitioned into matchings of $H(r,s,t+1)$ (the graph without
strips 1 and 4) and
$H(r+1,s+1,t-1)$ (the graph without strips 2 and 3), plus two additional
corner edges.  In both cases, it is possible to partition the collection
into matchings of $H(r+1,s+1,t)$
and $H(r,s,t)$.  Finally, it is impossible for the paths to run from
strips 1 to strip 3 and from strip 2 to strip 4 without intersecting.  Thus
\begin{eqnarray*}
\lefteqn{q^{(r+1)(s+1)s/2} P(r+1,s+1,t) \cdot q^{rs(s-1)/2}P(r,s,t) = } \\
& &      q^{s+t}\cdot q^{r(s+1)s/2} P(r,s+1,t)
  \cdot  q^{(r+1)s(s-1)/2}P(r+1,s,t) + \\
& &      q^{(r+1)(s+1)s/2} P(r+1,s+1,t-1) \cdot q^{rs(s-1)/2} P(r,s,t+1).
\end{eqnarray*}
Note how the factor of $q^{s+t}$ in the right-hand side comes from the edge of
weight $q^{s+t}$ that was not covered by either subgraph $H(r,s+1,t)$ or
$H(r+1,s,t)$.

We simplify this relation by dividing through by $q^{(r+1)(s+1)s/2+rs(s-1)/2}$
to get the desired relation:
\[P(r+1,s+1,t)P(r,s,t)
= q^t P(r,s+1,t) P(r+1,s,t) + P(r+1,s+1,t-1) P(r,s,t+1). \]
\end{proof}

Now we can prove MacMahon's formula for $P(a,b,c)$ by induction on $a+b+c$.  
When any of $a$, $b$, or $c$ are 0, $P(a,b,c)=1$.  Now suppose MacMahon's 
formula holds for all $a,b,c$ such that $a+b+c \leq r+s+t+1$.  We show 
MacMahon's formula holds for $(a,b,c)=(r+1,s+1,t)$:

\begin{eqnarray*}
P(r,s+1,t) P(r+1,s,t)
&=& \left( \prod_{i=1}^r \prod_{j=1}^s
\frac{1-q^{i+j+t-1}}{1-q^{i+j-1}} \right)^2
\left( \prod_{j=1}^s \frac{1-q^{j+r+t}}{1-q^{j+r}} \right)
\left( \prod_{i=1}^r \frac{1-q^{i+s+t}}{1-q^{i+s}} \right) \\
&=& P(r,s,t) \cdot \frac{1-q^{r+s-1}}{1-q^{r+s+t-1}}
\prod_{i=1}^{r+1} \prod_{j=1}^{s+1} \frac{1-q^{i+j+t-1}}{1-q^{i+j-1}} \\
P(r+1,s+1,t-1) P(r,s,t+1) &=& \frac{\prod_{i=1}^{r+1}
\prod_{j=1}^{s+1} (1-q^{i+j+t-2}) \prod_{i=1}^r \prod_{j=1}^s (1-q^{i+j+t})}
{\prod_{i=1}^{r+1} \prod_{j=1}^{s+1} (1-q^{i+j-1})\prod_{i=1}^r \prod_{j=1}^s
(1-q^{i+j-1})} \\
&=& \prod_{i=1}^r \prod_{j=1}^s \frac{1-q^{i+j+t-1}}{1-q^{i+j-1}}
\frac{(1-q^t)\prod_{i=1}^{r+1} \prod_{j=1}^s (1-q^{i+j+t-1}) \prod_{j=1}^r
(1-q^{j+s+t})}{\prod_{i=1}^r \prod_{j=1}^s 1-q^{i+j-1}} \\
&=& P(r,s,t) \cdot \frac{1-q^t}{1-q^{r+s+t-1}} \prod_{i=1}^{r+1}
\prod_{j=1}^{s+1} \frac{1-q^{i+j+t-1}}{1-q^{i+j-1}} \\
P(r+1,s+1,t) P(r,s,t) &=&
q^t P(r,s+1,t) P(r+1,s,t) + P(r+1,s+1,t-1) P(r,s,t+1) \\
&=& \frac{q^t(1-q^{r+s-1})+ (1-q^t)}{1-q^{r+s+t-1}} P(r,s,t)
\prod_{i=1}^{r+1} \prod_{j=1}^{s+1} \frac{1-q^{i+j+t-1}}{1-q^{i+j-1}} \\
&=& P(r,s,t)\prod_{i=1}^{r+1} \prod_{j=1}^{s+1}
\frac{1-q^{i+j+t-1}}{1-q^{i+j-1}}.
\end{eqnarray*}
Thus
\[ P(r+1,s+1,t) = \prod_{i=1}^{r+1} \prod_{j=1}^{s+1}
\frac{1-q^{i+j+t-1}}{1-q^{i+j-1}}. \]

\begin{figure}
\insertfig{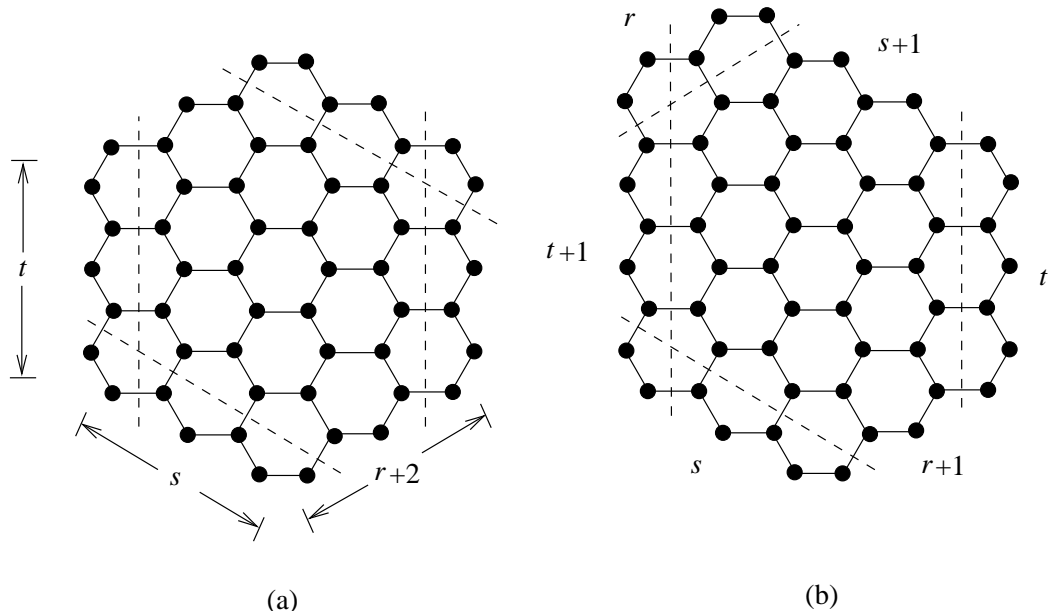}
\caption{Graphs for proving other bilinear relations.  (a) $r=1$, $s=t=3$.
(b) $r=2$, $s=3$, $t=3$.}
\label{fig:wtd-hexes}
\end{figure}

It is interesting to note a similarity between Theorems~\ref{thm:abcd}
and~\ref{thm:wtd-hex-1}.  In the proof of each theorem, the two paths always
run between vertices of opposite parity.
We can find additional bilinear relations with MacMahon's formula that are
analogous to Theorems~\ref{thm:ad-bc} and~\ref{thm:b-acd}.  For instance, if
we partition a hexagonal graph as shown in Figure~\ref{fig:wtd-hexes}(a), we
get
\begin{eqnarray*}
\lefteqn{q^{(r+2)s(s-1)/2} P(r+2,s,t) \cdot q^{rs(s-1)/2} P(r,s,t) = }\\
& &     (q^{(r+1)s(s-1)/2} P(r+1,s,t))^2 - \\
& &      q^{(r+1)(s-1)(s-2)/2}P(r+1,s-1,t+1) \cdot q^{(r+1)(s+1)s/2}
    P(r+1,s+1,t-1).
\end{eqnarray*}
After dividing through by $q^{(r+1)s(s-1)}$, we get
\[ P(r+2,s,t)P(r,s,t) = P(r+1,s,t)^2 - q^{r+1} P(r+1,s-1,t+1) P(r+1,s+1,t-1). \]

The relation analogous to Theorem~\ref{thm:b-acd} is:
\begin{eqnarray*}
\lefteqn{q^{rs(s-1)/2} P(r,s,t+1) \cdot q^{rs(s-1)/2} P(r,s,t) = } \\
& &      q^{(r+1)s(s-1)/2} P(r+1,s,t) \cdot q^{(r-1)s(s-1)/2} P(r-1,s,t+1) +\\
& &      q^{r(s+1)s/2} P(r,s+1,t) \cdot q^{r(s-1)(s-2)/2}P(r,s-1,t+1),
\end{eqnarray*}
which simplifies to
\[ P(r,s,t+1)P(r,s,t) = P(r+1,s,t)P(r-1,s,t+1)+q^r P(r,s+1,t) P(r,s-1,t+1). \]
Figure~\ref{fig:wtd-hexes}(b) shows how to prove this relation.  The graph
has sides $r,s+1,t,r+1,s,t+1$.  For each pair of hexagons, one is missing
one of the strips along the sides of length $r$, $t$, or $s$, and the other hexagon is missing
the other three strips (but contains the strip that the first hexagon
is missing).

By taking the limit as $q \rightarrow 1$, we derive relations
among the numbers of plane partitions fitting in $\bx{r,s,t}$.  These numbers
also enumerate rhombus tilings of semiregular hexagons with sides $r,s,t,r,s,t$.
In particular, the following relation was proven by
Doron Zeilberger in~\cite{z96}:
\[ N(r,s,t+1)N(r,s,t) = N(r+1,s,t)N(r-1,s,t)+N(r,s+1,t)N(r,s-1,t+1) \]
where $N(r,s,t)$ is $\lim_{q\rightarrow 1} P(r,s,t)$.

\section{Transpose Complement Plane Partitions}\label{sec:tcpp}
If we view a plane partition as a collection of stacks of cubes, certain
plane partitions will exhibit some symmetry.  Such symmetry classes are
outlined in~\cite{B99}.  The {\it complement} of a plane partition $\pi$ in
the box of dimensions $r \times s \times t$ is the set of cubes in the box
that are not in $\pi$, reflected through the center of the box.  A {\it
transpose complement plane partition (TCPP)\/} $\pi$ is one for which the
complement is the same as the reflection of $\pi$ in the plane $y=x$.  If we
visualize a TCPP as a rhombus tiling of a hexagon, the line of symmetry goes
through the midpoints of two sides of the hexagon.  Note that the sides of the
hexagon must be of the form $r, r, 2t, r, r, 2t$, and the line of symmetry goes
through the sides of length $2t$.  The following theorem about the number of
TCPPs was proved in~\cite{p88}:

\begin{theorem}\label{thm:tcpp}
The number of TCPPs in an $r \times r \times 2t$ box is
\[ {{t+r-1}\choose{r-1}} \prod_{1 \leq i \leq j \leq r-2}
\frac{2t+i+j+1}{i+j+1}. \]
\end{theorem}

Let $N(r,r,2t)$ be the the number of TCPPs in an $r \times r \times 2t$ box.

\begin{proposition}\label{prop:tcpp}
If $r \geq 2$ and $t \geq 1$, then
\[ N(r,r,2t)N(r-2,r-2,2t)=N(r-1,r-1,2t)^2 + N(r,r,2t-2)N(r-2,r-2,2t+2). \]
\end{proposition}

\begin{figure}
\insertfig{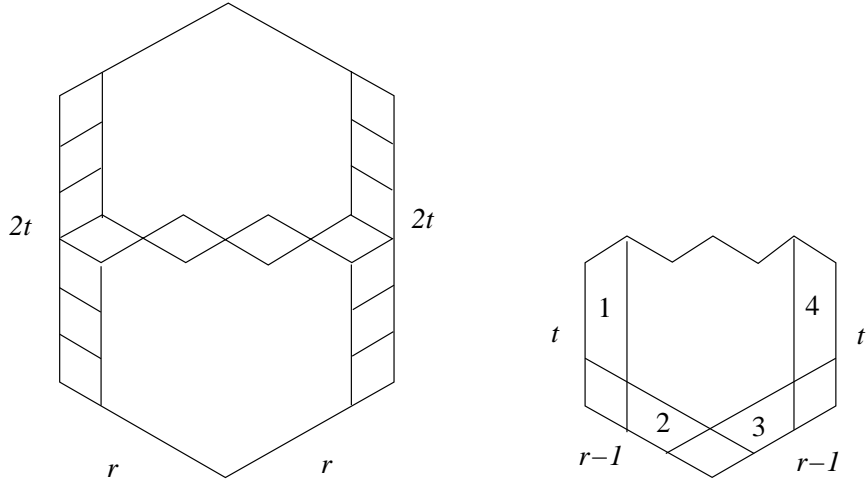}
\caption{Left: The forced rhombi in a TCPP.  Right: The four strips along
the sides of a (stripped) semihexagon.}
\label{fig:tcpp}
\end{figure}

\begin{proof}
Because of the symmetry of a TCPP, we only need to consider the number of
ways to tile one half of an $(r,r,2t,r,r,2t)$-hexagon.  Also note that the
triangles that lie on the line of symmetry must join to form rhombi.
We can cut the hexagon in half to form an $(r,t)$-semihexagon.  We can strip
this semihexagon even further since rhombi are forced along the sides of
length $t$. (This also shortens the sides of length $r$ by one.)  See
Figure ~\ref{fig:tcpp}, left.  Let us label the four strips of triangles 1, 2, 3,
and 4, so that strips 1 and 4 are along the sides of length $t$, and strips 2
and 3 are along the sides of length $r-1$ (see Figure~\ref{fig:tcpp}, right).
Removing all four strips would produce an $(r-2,t)$-semihexagon.
Removing only strips 1 and 2 (or only strips 3 and 4) would
produce a region with the same number of tilings as an $(r-1,t)$-semihexagon.
If we remove only strips 2 and 3, we shorten $t$ by one to form an $(r,t-1)$-
semihexagon (with its outer edges stripped).  Finally, if we remove only
strips 1 and 4, we get an $(r-1, t+1)$-semihexagon.  The relation follows
from graphical condensation.
\end{proof}

This relation was noted by Michael Somos in a private communication.

We can now prove Theorem~\ref{thm:tcpp} by induction on $r+t$.  The relevant
base cases are $N(r,r,0)=1$ for all $r$, and $N(0,0,2t)=1$ and $N(1,1,2t)=1$
for all $t$.  These cases are trivially established.  We now use
Proposition ~\ref{prop:tcpp} to prove the inductive step.  Given that the
formula holds for $N(r-2,r-2,2t), N(r-1,r-1,2t), N(r,r,2t-2)$, and $N(r-2,r-2,2t+2)$,
we show that it holds also for $N(r,r,2t)$.  We need to verify that
\begin{eqnarray*}
\lefteqn{{{t+r-1}\choose{r-1}}{{t+r-3}\choose{r-3}}
\prod_{1\leq i \leq j \leq r-2} \frac{2t+i+j+1}{i+j+1}
\prod_{1\leq i \leq j \leq r-4}\frac{2t+i+j+1}{i+j+1} =} \\
& & {{t+r-2}\choose{r-2}}^2 \left( \prod_{1\leq i \leq j \leq r-3}
\frac{2t+i+j+1}{i+j+1} \right)^2 + \\
& & {{t+r-2}\choose{r-1}}{{t+r-2}\choose{r-3}}\prod_{1\leq i \leq j \leq r-2}
\frac{2t-2+i+j+1}{i+j+1} \prod_{1\leq i \leq j \leq r-4}
\frac{2t+2+i+j+1}{i+j+1}.
\end{eqnarray*}
We divide the right hand side by the left hand side to obtain
\begin{eqnarray*}
\lefteqn{\frac{r-1}{t+r-1}\cdot\frac{t+r-2}{r-2}
\prod_{1 \leq i \leq r-2} \frac{i+(r-2)+1}{2t+i+(r-2)+1}
\prod_{1 \leq i \leq r-3} \frac{2t+i+(r-3)+1}{i+(r-3)+1}} \\
&+& \frac{t}{t+r-1}\cdot\frac{t+r-2}{t+1}
\prod_{1 \leq i \leq j \leq r-2} \frac{2t+i+j-1}{2t+i+j+1}
\prod_{1 \leq i \leq j \leq r-4} \frac{2t+i+j+3}{2t+i+j+1}.
\end{eqnarray*}
After a heavy dose of cancellations, this expression simplifies to
\begin{eqnarray*}
\lefteqn{ \frac{(r-1)(t+r-2)(2r-3)(2r-4)(2t+r-1)}{(t+r-1)(r-2)(r-1)(2t+2r-3)(2t+2r-4)} +
\frac{t(t+r-2)(2t+1)(2t+2)}{(t+r-1)(t+1)(2t+r-3)(2t+2r-4)} } \\
&=& \frac{(2r-3)(2t+r-1)}{(t+r-1)(2t+2r-3)}+\frac{t(2t+1)}{(t+r-1)(2t+2r-3)}
= \frac{(4rt+2r^2-5r-6t+3)+(2t^2+t)}{2t^2+2r^2+4rt-5t-5r+3} = 1
\end{eqnarray*}
and thus the inductive step and Theorem~\ref{thm:tcpp} follows.

\section{Acknowledgments}
Great thanks go to James Propp for suggesting the Holey Aztec Rectangle
problem, to which this paper owes its existence; for suggesting
some applications of graphical condensation, including domino tilings of 
fortresses and (non-Aztec) rectangles; and for editing and revising
this article, making many very helpful suggestions along the way.  
Thanks also to Henry Cohn for pointing out the similarity between the relations
in Theorems~\ref{thm:har} and Corollary~\ref{cor:4-cycle}, thus providing me 
an inspiration for graphical condensation; and for applying graphical 
condensation to 
placement probabilities in Aztec Diamonds.  Thanks to David Wilson
for creating the software that enabled the enumeration of tilings of Aztec 
rectangles, and finally thanks to the referees for providing some very helpful
comments for this article.

\end{document}